\documentclass[11pt,reqno,a4paper]{amsart}
\usepackage{euscript}

\newcommand{\N}{\ensuremath{\mathbb{N}}}

\newtheorem {theorem} {Theorem}

\newtheorem {lemma}  [theorem]{Lemma}

\begin{document}

\title[A note on the $3x+1$ conjecture]
{A note on the $3x+1$ conjecture}

\author[J. Llibre and C. Valls]
{Jaume Llibre$^1$ and Claudia Valls$^2$}

\address{$^1$ Departament de Matem\`{a}tiques, Universitat Aut\`{o}noma de Barcelona, 08193 Bellaterra, Barcelona, Catalonia, Spain}
\email{jllibre@mat.uab.cat}

\address{$^2$ Departamento de Matem\'atica, Instituto Superior T\'ecnico, Universidade de Lisboa, Av. Rovisco Pais 1049--001, Lisboa, Portugal} \email{cvalls@math.tecnico.ulisboa.pt}



\begin{abstract}
Let $g$ be a map from the set of positive integers into itself defined as follows: Let $x$ be a positive integer. If $x$ is odd, then $g(x)=3x+1$, and if $x$ is even, then $g(x)=x/2$. The $3x+1$ conjecture, also called the Collatz conjecture, states: For any positive integer $x$ there exists another positive integer $m$ such that the $m$-iterate of $x$ under the map $g$ is equal to $1$, i.e. $g^m(x)=1$. We provide some information related with this conjecture.
\end{abstract}

\maketitle

\section{Introduction and the main result}

The $3x+1$ or Collatz conjecture was stated by Collatz in the 1930's. For the history of the subject see the book of Lagarias \cite{La2} which describes methods (mainly statistical and computational) of approaching the conjecture and its generalizations, see also the book of Wirsching \cite{Wi} and the paper of Sinai \cite{Si} for an ergodic approach to the conjecture.

\begin{lemma}\label{L1}
The following statements hold.
\begin{itemize}
\item[(a)] If $p$ is even, then $2^p\equiv 1 \pmod 3$.
	
\item[(b)] If $p$ is odd, then $2^p\equiv 2 \pmod 3$.
\end{itemize}
\end{lemma}

We denote by $\N$ the set of positive integers (the set of natural numbers) and by $I$ the subset of $\N$ formed by the odd natural numbers. If $k\in I$, then $3k+1$ is even. So we can write $3k+1=2^n \bar k$, where $\bar k \in I$ and $n\in \N$. We define a function $f\colon I \to I$ by $f(k)=\bar k$, see \cite{La2}.

We tried to prove the conjecture using the induction method without success. We shall explain what we have done and where are the difficulties for a such kind of proof.

Let $E$ be the set of all $k\in I$ for which there exists $n\in \N$ such that $f^n(k) = 1$. Here $f^n(k)$ denotes the $n$th-iterate of $k$ under the function $f$. To prove the $3x+1$ conjecture is equivalent to prove that $E=I$. It is known that every $k\in I$ with $k<2^{60}$ belongs to $E$, see \cite{Ba}.

Our intent to find a proof for induction was as follows. For each $k\in I$ with $k\ge 2^{60}$ we assume that all odd integers less than or equal to $ k-2$ are in $E$, and we want to prove that $k\in E$.

Let $k=2^p h-1$ with $p\in \N$ and $h\in I$ be such an integer.

\noindent{\it Case} 1: $p=1$. Since $k=2h-1>1$, we have $3k+1=2(3h-1)=2^q \bar h$ with $q\ge 2$ and $\bar h$ not divisible by $2$. Therefore,
$
f(k)=f(2h-1)= \bar h\le \frac{3h-1}{2}< 2h-1=k.
$
Since $f(k)<k$, by the induction hypothesis $f(k)\in E$. So there exists $m\in \N$ such that $f^m(f(k))=1$. Hence $f^{m+1}(k)=1$ and $k\in E$.

\smallskip

\noindent{\it Case} 2: $p\ge 2$ and $h\equiv 0 \pmod 3 $. Then $h=3\bar h$ for some $\bar h\in I$ and so $k=2^p h-1= 2^p 3 \bar h-1$. Let $\bar k= 2^{p+1}\bar h-1$.  Then $3\bar k+1= 2(2^p h-1)$. Therefore, $f(\bar k)= 2^p h-1= k$. Since $\bar k<k$, by the induction hypothesis $\bar k\in E$ and so there exists $m\in \N$ such that $f^m(\bar k)=1$. Hence $f^{m-1}(k)=1$~and $k\in E$.	

\smallskip

\noindent{\it Case} 3: $p\ge 2$, $h\equiv 1 \pmod 3$ and $p$ odd. Let $h=3\ell+1$. We claim that 
$$
m=\frac{2^2k-1}{3}=\frac{2^2(2^p(3\ell+1)-1)-1}{3}
$$
is an odd integer such that $f(m)=k$. Clearly that the claim will be proved if $m$ is an integer. From the definition of $m$ we can write
$$
3m+1=2^2(2^p(3\ell+1)-1)=2^{p+2}\cdot 3\ell+2^2(2^p-1).
$$
Taking in this equality modulo $3$ we obtain $1\equiv 2^2(2^p-1)$. Since $p$ is odd from Lemma \ref{L1}(b) we get that  $1\equiv 2^2$. Since this last equivalence holds we have that $m$ is an integer. So the claim is proved.

Now we claim that
$$
r=\frac{4m-1}{3}
$$
is an odd integer such that $f(r)=m$. Again this claim will be proved if $r$ is an integer. From the definition of $r$ we can write
$$
4m= 2^2\frac{2^2(2^p(3\ell+1)-1)-1}{3}=3r+1,
$$
or equivalently
$$
2^{p+4}\cdot3\ell+2^{p+4}-2^4-2^2=9r+3.
$$
Taking in this equality modulo $3$ we obtain $2^{p+4}-20\equiv 0$, i.e.
$2^{p+4}-2\equiv 0$, and this last equivalence holds by Lemma \ref{L1}(b) because $p$ is odd. Therefore $r$ is an integer.

We have that
$$
k=\frac{3m+1}{4}=\frac{3\frac{3r+1}{4}+1}{4}=\frac{9r+7}{16}>\frac{9r}{16}>r.
$$
Since $f^2(r)=f(m)=k$ and $r<k$ by induction assumption $r\in E$. So there exists $n$ such that $f^n(r)=1$. Therefore $f^{n-2}(k)=1$ and $k\in E$.

\smallskip

\noindent{\it Case} 4: $p\ge 2$, $h\equiv 2 \pmod 3$ and $p$ even. Let $h=3\ell+2$. We claim that 
$$
m=\frac{2^2k-1}{3}=\frac{2^2(2^p(3\ell+2)-1)-1}{3}
$$
is an odd integer such that $f(m)=k$. Clearly that the claim will be proved if $m$ is an integer. From the definition of $m$ we can write
$$
3m+1=2^2(2^p(3\ell+2)-1)=2^{p+2}\cdot 3\ell+2^2(2^{p+1}-1).
$$
Taking in this equality modulo $3$ we obtain $1\equiv 2^2(2^{p+1}-1)$. Since $p$ is even from Lemma \ref{L1}(b) we get that  $1\equiv 2^2$. Since this last equivalence holds we have that $m$ is an integer. So the claim is proved.

Now we claim that
$$
r=\frac{4m-1}{3}
$$
is an odd integer such that $f(r)=m$. Again this claim will be proved if $r$ is an integer. From the definition of $r$ we can write
$$
4m= 2^2\frac{2^2(2^p(3\ell+2)-1)-1}{3}=3r+1,
$$
or equivalently
$$
2^{p+4}\cdot3\ell+2^{p+5}-2^4-2^2=9r+3.
$$
Taking in this equality modulo $3$ we obtain $2^{p+5}-20\equiv 0$, i.e.
$2^{p+5}-2\equiv 0$, and this last equivalence holds by Lemma \ref{L1}(b) because $p$ is even. Therefore $r$ is an integer.

Again we have that
$$
k=\frac{3m+1}{4}=\frac{3\frac{3r+1}{4}+1}{4}=\frac{9r+7}{16}>\frac{9r}{16}>r.
$$
Since $f^2(r)=f(m)=k$ and $r<k$ by induction assumption $r\in E$. So there exists $n$ such that $f^n(r)=1$. Therefore $f^{n-2}(k)=1$ and $k\in E$.

\smallskip

\smallskip

\noindent{\it Case} 5: $p\ge 2$, $h\equiv 1 \pmod 3$ and $p$ even. So we can write $h=3\ell+1$. We look for an
$$
m=\frac{2^s k-1}{3}=\frac{2^s(2^p(3\ell+1)-1)-1}{3}
$$
integer such that $f(m)=k$ with $s\in \N$. Clearly that if $m$ exists then it is an odd integer. From the definition of $m$ we can write
$$
3m+1=2^s(2^p(3\ell+1)-1)=2^{p+s}\cdot 3\ell+2^s(2^p-1).
$$
Taking in this equality modulo $3$ since $p$ is even we obtain that $1\equiv 0$, a contradiction. So such $m$ does not exist. Hence for the integers $k$ in case 5 there does not exist an integer $\bar k<k$ such that $f(\bar k)= k$.

\smallskip

\noindent{\it Case} 6: $p\ge 2$, $h\equiv 2 \pmod 3$ and $p$ odd. So we can write $h=3\ell+2$.  We again look for an
$$
m=\frac{2^sk-1}{3}=\frac{2^s(2^p(3\ell+2)-1)-1}{3}
$$
integer such that $f(m)=k$. Clearly that if $m$ exists then it is an integer. From the definition of $m$ we can write
$$
3m+1=2^s(2^p(3\ell+2)-1)=2^{p+s}\cdot 3\ell+2^s(2^{p+1}-1).
$$
Taking in this equality modulo $3$  since $p$ is even from Lemma \ref{L1}(b) we get that  $1\equiv 0$, a contradiction. So such $m$ does not exist. Hence for the integers $k$ in case 6 there does not exist an integer $\bar k<k$ such that $f(\bar k)= k$. 

\smallskip

If the odd integer $k$ belongs to the cases 5 or 6, then
$$
f^n(k)=f^n(2^ph-1)=3^n2^{p-n}h-1,\ \mbox{for $n=1,\ldots,p-1$.}
$$
So
$$
k<f(k)<f^2(k)<\ldots f^{p-1}(k), \ \mbox{and $f^p(k)=3^ph-1$.}
$$

Here appears the problem in order to try to end the proof by induction, how to control the iterates under $f$ of the integers of the form $3^ph-1$ with $p\ge 2$ and $h$ odd in the cases 5 and 6.

\end{document}